\documentclass[psamsfonts, 12pt, leqno]{amsart}

\usepackage{hyperref}

\newtheorem{Theorem}{Theorem}
\newtheorem{Corollary}[Theorem]{Corollary}
\newtheorem{Lemma}[Theorem]{Lemma}

\def\eps{\varepsilon}
\def\intl{\int\limits}

\begin{document}
\title[Self-approximation of Dirichlet $L$-functions ]
{Self-approximation of Dirichlet $L$-functions}

\author{\sc Ram\={u}nas Garunk\v{s}tis}

\address{Ram\={u}nas Garunk\v{s}tis \\
Department of Mathematics and Informatics, Vilnius University \\
Naugarduko 24, 03225 Vilnius, Lithuania}

\thanks{Supported by grant No MIP-94 from the Research Council of Lithuania}

\email{ramunas.garunkstis@mif.vu.lt}

\urladdr{www.mif.vu.lt/~garunkstis}

\date{}

\setcounter{equation}{0}


\begin{abstract}
Let $d$ be a real number, let $s$ be in a fixed
compact set of the strip $1/2<\sigma<1$, and let
$L(s, \chi)$ be the Dirichlet $L$-function. The
hypothesis  is that for any real number $d$ there
exist `many' real numbers $\tau$ such that  the
shifts
 $L(s+i\tau, \chi)$ and $L(s+id\tau, \chi)$ are `near' each other. If $d$ is an algebraic irrational number then this was obtained by T. Nakamura.  \L. Pa\'nkowski solved the case then $d$ is a transcendental number. We prove the case then $d\ne0$ is a rational number.  If $d=0$ then by B. Bagchi we know that the above hypothesis is equivalent to the Riemann hypothesis for the given Dirichlet $L$-function. We also consider a more general  version of the above problem.
\end{abstract}

\maketitle

\section{Introduction}

Let, as usual, $s=\sigma+it$ denote a complex
variable. For $\sigma>1$, the Dirichlet
$L$-function is given by
\[
L(s,\chi)=\sum_{n=1}^{\infty}\frac{\chi(n)}{n^s},
\]
where $\chi(n)$ is a Dirichlet character $\bmod\,
q$.  For $q=1$ we get $L(s,\chi)=\zeta(s)$, where
$\zeta(s)$ is the Riemann zeta-function.

In \cite{bohr} Bohr proved that if $\chi$ is a
nonprincipal character, then the Riemann
hypothesis for $L(s, \chi)$ is equivalent to the
almost periodicity of $L(s, \chi)$ in the half
plane $\sigma>1/2$. A function $f(s)$ is  {\it
almost periodic} in a region $E\subset\Bbb C$ if
for any positive $\eps$ and any compact subset
$K$ in $E$ there exists a sequence of real
numbers $\dots<\tau_{-1}<0<\tau_1<\tau_2<\dots$
such that
$$\liminf_{m\to\pm\infty}(\tau_{m+1}-\tau_m)>0,\qquad\limsup_{m\to\pm\infty}\frac{\tau_m}m<\infty$$
and
$$|f(s+i\tau_m)-f(s)|<\eps\quad\textup{for all}\ s\in K \ \textup{and}\ m\in\Bbb Z$$
 hold.
Bohr \cite{bohr} also obtained that every
Dirichlet series is almost-periodic in its
half-plane of absolute convergence. Effective
upper bounds for the almost periodicity of
Dirichlet series with Euler products in the
half-plane of absolute convergence were
considered by Girondo and Steuding \cite{gs}.
Note that  every Dirichlet $L$-function is almost
periodic in the sense of Besicovitch on any
vertical line of the strip $1/2<\sigma<1$. For
this and related results see Besicovitch
\cite{bes} and Mauclaire \cite{mauc},
\cite{mauc2}.

Bagchi \cite{bag} proved that the Riemann
hypothesis for $L(s, \chi)$ ($\chi$  is an
arbitrary Dirichlet character) is true if and
only if for any compact subset $\mathcal K$ of
the strip $1/2<\sigma<1$ and for any $\eps>0$
\begin{align}\label{strong}
\liminf_{T\to\infty}\frac1
T\text{meas}\left\{\tau\in[0, T] :
\max_{s\in\mathcal K}|L(s+i\tau, \chi)-L(s,
\chi)|<\eps\right\}>0,
\end{align}
where $\text{meas}\, A$ stands for the Lebesgue
measure of a measurable set $A$. Bagchi says that
the Dirichlet $L$-function $L(s, \chi)$ is {\it
strongly recurrent} on the strip
$\sigma_0<\sigma<\sigma_1$ if  (\ref{strong}) is
valid for any compact $\mathcal K$ of the strip
$\sigma_0<\sigma<\sigma_1$. The strong recurrence
is connected with the universality property of
Dirichlet series. More about the universality and
the strong recurrence see Bagchi  \cite{bagdis},
\cite{bag}, \cite{bag87}, and Steuding
\cite{steu}.

There are several unconditional results
concerning the self-approximation of Dirichlet
$L$-functions in the critical strip. Let
$\mathcal K$ be a compact subset of the strip
$1/2<\sigma<1$ and let $\lambda\in\Bbb R$ be such
that $\mathcal K$ and  $\mathcal K+i\lambda
:=\{s+i\lambda : s\in\mathcal K\}$ are disjoint.
From Kaczorowski, Laurin\v cikas and Steuding
\cite{kls} it follows that for any character
$\chi$ and any $\eps>0$
\begin{align*}
\liminf_{T\to\infty}\frac1
T\text{meas}\left\{\tau\in[0, T] :
\max_{s\in\mathcal K}|L(s+i\lambda+i\tau,
\chi)-L(s+i\tau, \chi)|<\eps\right\}>0.
\end{align*}
Nakamura \cite{nak} considered the joint
universality of shifted Dirichlet $L$-functions.
His Theorem 1.1 leads to the following statement.
If $1=d_1, d_2, \dots, d_m$ are algebraic real
numbers {\it linearly independent} over $\Bbb Q$,
then for any Dirichlet character $\chi$ and any
$\eps>0$
\begin{align}\label{naka}
\liminf_{T\to\infty}&\frac1
T\text{meas}\left\{\tau\in[0, T] \right.
:\\&\left. \max_{1\le j,k\le m}\max_{s\in\mathcal
K}|L(s+id_j\tau, \chi)-L(s+id_k\tau,
\chi)|<\eps\right\}>0.\nonumber
\end{align}
If $m=2$ then Pa\'nkowski \cite{pan} using Six
Exponentials Theorem showed that (\ref{naka})
holds for $d_1, d_2$ are real numbers linearly
independent over $\Bbb Q$.

We prove the following theorem.
\begin{Theorem}\label{th}
Let $1=d_1, d_2, \dots, d_m$ be nonzero algebraic
real numbers and let $\mathcal K$ be a compact
subset of the strip $1/2<\sigma<1$. Then for any
Dirichlet character $\chi$ and any $\eps>0$ the
inequality (\ref{naka}) is valid.
\end{Theorem}
Note that Theorem \ref{th} remains true if $d_1,
d_2, \dots, d_m$ are replaced by $dd_1, dd_2,
\dots, dd_m$, where $d\in\Bbb R$. The next
theorem shows that `$\liminf $' in the inequality
(\ref{naka})  often can be replaced by `$\lim$'.
\begin{Theorem}\label{th2}
Let $d_1, d_2, \dots, d_m$ be any real numbers,
let $\chi_1, \chi_2, \dots,\chi_m$ be  any
Dirichlet characters, and  let $\mathcal K$ be a
compact subset of the strip $1/2<\sigma<1$. Then
for any $\eps>0$, except an at most countable set
of $\eps$, there exists a limit
\begin{align*}
\lim_{T\to\infty}&\frac1
T\text{meas}\left\{\tau\in[0, T] : \max_{1\le
j,k\le m}\max_{s\in\mathcal K}|L(s+id_j\tau,
\chi_j)-L(s+id_k\tau,
\chi_k)|<\eps\right\}.\nonumber
\end{align*}
\end{Theorem}

The mentioned results of Nakamura and Pa\'nkowski
together with Theorem \ref{th} and Theorem
\ref{th2} lead to the following corollary.

\begin{Corollary}\label{cor}
Let $d$ be a nonzero real number and let
$\mathcal K$ be a compact subset of the strip
$1/2<\sigma<1$. Then for any Dirichlet character
$\chi$ and any $\eps>0$, except an at most
countable set of $\eps$,
\begin{align}\label{lim}
\lim_{T\to\infty}\frac1
T\text{meas}\left\{\tau\in[0, T] :
\max_{s\in\mathcal K}|L(s+i\tau,
\chi)-L(s+id\tau, \chi)|<\eps\right\}>0.
\end{align}
\end{Corollary}

From the proof of Theorem \ref{th2} we see that
for any real numbers $d_1,\dots,d_m$ and for any
Dirichlet characters $\chi_1,\dots,\chi_m$ the
function
\[
 g(\tau)=
\max_{1\le j,k\le m}\max_{s\in\mathcal
K}\left|L(s+id_j\tau, \chi_j)-L(s+id_k\tau,
\chi_k) \right|
\]
is Besicovitch almost periodic function (for the
definition see Section \ref{secproof} above the
proof of Theorem \ref{th2}). Let $\eps>0$ be such
that the limit (\ref{lim}) exists. For such $\eps$ we
define a characteristic function
\(I_\eps(\tau)\), $\tau\in\Bbb R$, by
\begin{align}
I_\eps(\tau)=\left\{
\begin{matrix}
1,&\text{if}&g(\tau)<\eps,\\
0,&\text{if}&g(\tau)\ge\eps.
\end{matrix}
\right.
\end{align}
It is known (Jessen and A. Wintner \cite[Section
12]{jw})
 that \(I_\eps(\tau)\) is
Besicovitch almost periodic function also. Thus
we can say that self-approximations of Dirichlet
$L$-functions, considered in this paper,
 usually appear in a regular way.

 Theorem \ref{th} and Theorem \ref{th2} are
proved in Section \ref{secproof}. Next we state
several lemmas.

\section{Lemmas}

We start from the following statement.
\begin{Lemma}\label{03}
Let $\mathcal K$ be a compact subset of the
rectangle $U$. Let
$$d=\min_{z\in\partial\, U}\min_{s\in \mathcal K}|s-z|.$$
If $f(s)$ is analytic on $U$ and
$$\int_{U}\left|f(s)\right|^2 d\sigma dt\le\eps,$$
then
$$\max_{s\in \mathcal K}\left|f(s)\right|\le\frac{\sqrt{\eps/\pi}}d.$$
\end{Lemma}
\proof The lemma can be found in Gonek \cite{gon}
(Lemma 2.5).
\endproof

\begin{Lemma}\label{13}
Let $a_1,\dots,a_N$ be real numbers linearly
independent over the rational numbers. Let
$\gamma$ be a region of the $N$-dimensional unit
cube with volume $V$ (in the Jordan sense). Let
$I_\gamma(T)$ be the sum of the intervals between
$t=0$ and $t=T$ for which the point
$(a_1t,\dots,a_Nt)$ is$\mod 1$ inside $\gamma$.
Then
$$\lim_{T\to\infty}\frac{I_\gamma(T)}{T}=V.$$
\end{Lemma}
\proof This is Theorem 1 in Apendix, Section 8,
of Voronin and Karatsuba \cite{vor}.
\endproof

For a curve $\omega(t)$ in $\Bbb R^N$ we
introduce the notation
$$\left\{\omega(t)\right\}=\left(\omega_1(t)-[\omega_1(t)],\dots,\omega_N(t)-[\omega_N(t)]\right),$$
where $[x]$ denotes the integral part of
$x\in\Bbb R$.
\begin{Lemma}\label{04}
Suppose that the curve $\omega(t)$ is uniformly
distributed$\mod 1$ in $\Bbb R^N$. Let $D$ be a
closed and Jordan measurable subregion of the
unit cube in $\Bbb R^N$and let $\Omega$ be a
family of complex-valued continuos functions
defined on $D$. If $\Omega$ is uniformly bounded
and equicontinuous, then
$$\lim_{T\to\infty}\frac{ 1}{T}\int_0^T f\left(\{\omega(t)\}\right){ 1}_D(t)dt=\int_D f(x_1,\dots,x_N) d x_1\dots d x_N$$
uniformly with respect to $f\in\Omega$, where
$1_D(t)$ is equal to $1$ if $\omega(t)\in D \mod
1$, and $0$ otherwise.
\end{Lemma}
\proof The lemma is Theorem 3 in  Appendix,
Section 8, of Voronin and Karatsuba \cite{vor}.
\endproof

\begin{Lemma}\label{indep}
Let $p_n$ be the $n$th prime number and $1=d_1,
d_2,\dots,d_l$ be algebraic real numbers which
are linearly independent over $\Bbb Q$. Then the
set $\{d_k\log p_n\}_{n\in\Bbb N}^{1\le k\le l}$
is linearly independent over $\Bbb Q$.
\end{Lemma}
\proof This is Proposition 2.2 in Nakamura
\cite{nak}. The proof is based on Baker's
\cite[Theorem 2.4]{bak} result.
\endproof

\section{Proof of  Theorem \ref{th} and Theorem \ref{th2}}\label{secproof}

\proof[Proof of Theorem \ref{th}] We define a
truncated Dirichlet $L$-function
$$L_v(s, \chi)=\prod_{p\le v}\left(1-\frac{\chi(p)}{p^s}\right)^{-1}.$$
Roughly speaking, we first prove Theorem \ref{th}
for the truncated Dirichlet $L$-function and
later we show that the tail is small.

 Let $\{d_1,
d_2,\dots, d_l\}$ be a maximal linearly
independent (over $\Bbb Q$) subset of the set
$\{d_1, d_2,\dots, d_m\}$. Then there are
integers $a\ne0$ and $a_{k,1}, a_{k,2},\dots,
a_{k,l}$  such that
\begin{align}\label{lin}
 d_k=\frac1 a\left(a_{k,1}d_1+a_{k,2}d_2+\dots+a_{k,l}d_l\right) \quad\textup{for}\quad l< k\le m.
\end{align}

Let
$$A=\max_{{l<k\le m }}\{|a_{k, 1}|+|a_{k, 2}|+\dots+|a_{k, l}|\}.$$
 Denote by
$\Vert x\Vert$ the minimal distance of $x\in\Bbb
R$ to an integer. If
\begin{align}\label{tau}
\left\Vert\tau\frac{d_n\log p}{2\pi
a}\right\Vert<\delta \quad\textup{for}\quad p\le
v\ \textup {and}\ 1\le n\le l
\end{align}
then
\begin{align*}
\left\Vert\tau\frac{d_n\log p}{2\pi
}\right\Vert<a\delta \quad\textup{for}\quad p\le
v\ \textup {and}\ 1\le n\le l
\end{align*}
and, by the relation (\ref{lin}),
$$\left\Vert \tau\frac{d_k\log p}{2\pi}\right\Vert<A\delta \quad\textup{for}\quad p\le v\ \textup {and}\ l< k\le m.$$
 By this and by
the continuity of the function $L_v(s, \chi)$ we
have that for any $\eps>0$ there is $\delta>0$
such that for $\tau$ satisfying (\ref{tau})
 \begin{align}\label{maxmax}
 \max_{1\le k,n\le m}\max_{s\in\mathcal K}
 \left|\log L_v(s+id_k\tau, \chi)
 -\log L_v(s+id_n\tau, \chi)\right|<\eps.
 \end{align}

For  positive numbers $\delta$, $v$, and $T$ we
define the set
\begin{align}\label{st}
S_T=S_T(\delta, v)=\left\{\tau : \tau\in[0,T],
\left\|\tau\frac{d_n\log p}{2\pi
a}\right\|<\delta, \
 p\le v, \ 1\le n \le l\right\}.
\end{align}
Let $U$ be an open bounded rectangle with
vertices on the lines $\sigma=\sigma_1$ and
$\sigma=\sigma_2$, where
$1/2<\sigma_1<\sigma_2<1$, such that the set
$\mathcal K$ is in $U$. Let $y>v$. We have
\begin{align*}
&\frac1T\intl_{S_T}\intl_{U}\sum_{k=1}^m\left|
\log L_y(s+id_k\tau, \chi)- \log L_v(s+id_k\tau,
\chi)\right|^2 d\sigma dt d\tau
\\ &=\sum_{k=1}^m\intl_{U}\frac1T\intl_{S_T}\left|\log L_y(s+id_k\tau, \chi)- \log L_v(s+id_k\tau, \chi)\right|^2 d\tau d\sigma
dt.
\end{align*}
For the inner integrals of the right-hand side of
the last equality  we will apply Lemma \ref{04}.
Let $p_n$ be the $n$th prime number. There are
indexes $M$ and $N$ such that $p_M\le v<p_{M+1}$
and   $p_N\le y<p_{N+1}$. By  generalized
Kronecker's theorem (Lemma \ref{13}) and by Lemma
\ref{indep} the curve
$$\omega(\tau)=\left(\tau \frac{d_k\log p_n}{2\pi a}\right)_{1\le n\le N}^{1\le k\le l}$$
is uniformly distributed$\mod1$ in $\Bbb R^{lN}$.
Let $R'$ be a subregion of the $lN$-dimensional
unit cube defined by inequalities
$$\|y_{k,n}\|\le\delta\quad\text{for}\quad 1\le k\le l\ \text{and}\ 1\le n\le M$$
 and
$$\left|y_{k,n}-\frac12\right|\le\frac12 \quad\text{for}\quad 1\le k\le l\ \text{and}\ M+1\le n\le N.$$
 Let $R$ be a subregion of the $lM$-dimensional unit cube defined by inequalities
$$\|y_{k,n}\|\le\delta\quad\text{for}\quad 1\le k\le l\ \text{and}\ 1\le n\le M.$$
Clearly
$$\operatorname{meas}R'=\operatorname{meas}R=(2\delta)^{lM}.$$
Note that
\begin{align}\label{logL}
&\log L_y(s+id_k\tau, \chi)- \log L_v(s+id_k\tau,
\chi)=\log\frac{L_y}{L_v}(s+id_k\tau,
\chi)\nonumber
\\
&=-\sum_{v<p\le
y}\log\left(1-\frac{\chi(p)}{p^{s+id_k\tau}}\right)=\sum_{v<p\le
y}\sum_{j=1}^\infty\frac{\chi^j(p)}{jp^{j(s+id_k\tau)}}
\\
&=\sum_{M<n\le
N}\sum_{j=1}^\infty\frac{\chi^j(p_n)}{jp_n^{j(s+id_k\tau)}}.\nonumber
\end{align}
Thus in view of the linear dependence (\ref{lin})
we get
\begin{align*}
&\lim_{T\to\infty}
\frac1T\int_{S_T}\sum_{k=1}^m\left|\log\frac{L_y}{L_v}(s+id_k\tau,
\chi)\right|^2 d\tau
\\
& =
\lim_{T\to\infty}\frac1T\int_{S_T}\left(\sum_{k=1}^l\left|\log\frac{L_y}{L_v}(s+id_k\tau,
\chi)\right|^2\right.
\\
&
\left.\phantom{=}+\sum_{k={l+1}}^m\left|\log\frac{L_y}{L_v}\left(s+\frac{i}{a}(a_{k,1}d_1+a_{k,2}d_2+\dots+a_{k,l}d_l)\tau,
\chi\right)\right|^2\right) d\tau.
\end{align*}
By Lemma \ref{04} and equality (\ref{logL}) we
obtain that the last limit is equal to
\begin{align*}
& \int_{R'}\left(\sum_{k=1}^l\left|\sum_{M<n\le
N}\sum_{j=1}^\infty\frac{\chi^j(p)e^{-2\pi
ijay_{k,n}}}{jp_n^{js}}\right|^2 \right.
\\&
\left.+\sum_{k=l+1}^m\left|\sum_{M<n\le
N}\sum_{j=1}^\infty\frac{\chi^j(p)e^{2\pi
ij(a_{k,1}y_{1,n}+a_{k,2}y_{2,n}+\dots+a_{k,l}y_{l,n})}}{jp_n^{js}}\right|^2\right)dy_{1,1}\dots
dy_{l,N}
\\
& =\operatorname{meas}
R\intl_0^1\dots\intl_0^1\left(\sum_{k=1}^l\left|\sum_{M<n\le
N}\sum_{j=1}^\infty\frac{\chi^j(p)e^{-2\pi
ijay_{k,n}}}{jp_n^{js}}\right|^2 \right.
\\&
\left.+\sum_{k=l+1}^m\left|\sum_{M<n\le
N}\sum_{j=1}^\infty\frac{\chi^j(p)e^{2\pi
ij(a_{k,1}y_{1,n}+a_{k,2}y_{2,n}+\dots+a_{k,l}y_{l,n})}}{jp_n^{js}}\right|^2\right)dy_{1,M+1}\dots
dy_{l,N}
\\
&= m\operatorname{meas} R \sum_{v<p\le
y}\sum_{j=1}^\infty\frac{1}{jp^{2j\sigma}}\ll
\operatorname{meas} R
\sum_{p>v}\frac{1}{p^{2\sigma}}.
\end{align*}
Consequently
\begin{align}\label{iint}
&\frac1T\intl_{S_T}\intl_{U}\sum_{k=1}^m\left|
\log L_y(s+id_k\tau, \chi)- \log L_v(s+id_k\tau,
\chi)\right|^2 d\sigma dt d\tau
\\
& \ll \operatorname{meas} R
\sum_{p>v}\frac{1}{p^{2\sigma_1}}.\nonumber
\end{align}
 Again, by generalized
Kronecker's theorem (Lemma \ref{13}),
\begin{align}\label{kron}
\lim_{T\to\infty}\frac1T\operatorname{meas}
S_T=\operatorname{meas} R.
\end{align}
By (\ref{iint}) and (\ref{kron}), for large $v$,
as $T\to\infty$, we have
\begin{align*}
\operatorname{meas}\left\{\tau : \tau\in S_T,
\int_{
U}\sum_{k=1}^m\left|\log\frac{L_y}{L_v}(s+id_k\tau,
\chi)\right|^2 d\sigma
dt<\sqrt{\sum_{p>v}\frac{1}{p^{2\sigma_1}}}\right\}
>\frac12T\operatorname{meas}R.
\end{align*}
Then Lemma \ref{03} gives
\begin{align*}
&\operatorname{meas}\left\{\tau : \tau\in S_T,
\max_{s\in\mathcal
K}\sum_{k=1}^m\left|\log\frac{L_y}{L_v}(s+id_k\tau,
\chi)\right|^2
d\tau\le\frac1{d\sqrt{\pi}}\left(\sum_{p>v}\frac{1}{p^{2\sigma_1}}\right)^\frac14\right\}
\\
&>\frac12T\operatorname{meas}R,
\end{align*}
where $d=\min_{z\in\partial\, U}\min_{s\in
\mathcal K}|s-z|.$ By the continuity of the
logarithm we obtain that for  any $\eps>0$ there
is $v=v(\eps)$ such that for any $y>v$
\begin{align}\label{maxsum}
&\operatorname{meas}\left\{\tau : \tau\in S_T,
\max_{s\in\mathcal
K}\sum_{k=1}^m\left|L_y(s+id_k\tau,
\chi)-L_v(s+id_k\tau, \chi)\right|^2
d\tau<\eps\right\}
\\
&>\frac12T\operatorname{meas}R.\nonumber
\end{align}
Now we will prove that for any $\delta>0$ there
is $y=y(\delta)$ such that
\begin{align}\label{maxsum2}
&\operatorname{meas}\left\{\tau : \tau\in [0,T],
\max_{s\in\mathcal
K}\sum_{k=1}^m\left|L(s+id_k\tau,
\chi)-L_y(s+id_k\tau, \chi)\right|^2
d\tau<\delta\right\}
\\
&>(1-\delta)T.\nonumber
\end{align}
The last formula together with (\ref{maxmax}),
(\ref{st}) and (\ref{maxsum}) yields
Theorem~\ref{th}. We return to the proof of
(\ref{maxsum2}). By the mean value theorem of the
Dirichlet $L$-function (Steuding \cite{steu},
Corollary 6.11) and by Carlson's Theorem
(Titchmarsh \cite{titch}, Chapter 9.51) we obtain
\begin{align*}
\lim_{T\to\infty}\frac1T\intl_{0}^{T}\left|L(s+ix\tau,
\chi)-L_{y}(s+ix\tau, \chi)\right|^2
 d\tau=\sum_{n>y}\frac{|\chi(n)|}{n^{2\sigma}},
\end{align*}
where $x$ is fixed. Thus (\ref{maxsum2}) follows
in view of
\begin{align*}
\intl_{0}^{T}\intl_{U}\sum_{k=1}^m\left|L(s+id_k\tau,
\chi)-L_{y}(s+id_k\tau, \chi)\right|^2 d\sigma dt
d\tau\ll
\sum_{n>y}\frac{|\chi(n)|}{n^{2\sigma_1}}.
\end{align*}
Theorem \ref{th} is proved.
\endproof
The proof of Theorem \ref{th2} is based on the
ideas of Mauclaire \cite{mauc}, \cite{mauc2}. It
uses the theory of Besicovitch almost periodic
functions. We recall related definitions.

Let
$$P(\tau)=\sum_{n\in F} a_ne^{i\lambda_n \tau},$$
where $F$ is a finite set, $\lambda_n$ are any
real numbers, and the coefficients $a_n$ are any
complex numbers. For real  $\tau$ we say that
$P(\tau)$ is a {\it trigonometric polynomial}.

A  function $f : \Bbb R\to\Bbb C$ is called {\it
uniformly almost periodic} ($U.A.P.$) if given
any $\eps>0$, there exists a trigonometric
polynomial $P(\tau)$ such that
\[\sup_{\tau\in\Bbb R}|f(\tau)-P(\tau)|\le\eps.\]

A  function $f : \Bbb R\to\Bbb C$ is called {\it
$B^q$  almost periodic} ($B^q.A.P.$), $q\ge1$, if
given any $\eps>0$, there exists a trigonometric
polynomial $P(\tau)$ such that
\begin{align}\label{bap}
\limsup_{T\to\infty}\frac1{2T}
\int_{-T}^{T}|f(\tau)-P(\tau)|^q d\tau\le\eps.
\end{align}
If $q=1$ then we write $B.A.P.$ (Besikovitch
almost periodic) instead of $B^1.A.P.$ For any
$q\ge1$ it is clear that every $U.A.P.$ function
is $B^q.A.P.$ and that every $B^q.A.P.$ function
is $B.A.P.$

\proof[Proof of Theorem \ref{th2}] Let
\[
 g(\tau)=
\max_{1\le j,k\le m}\max_{s\in\mathcal
K}\left|L(s+id_j\tau, \chi_j)-L(s+id_k\tau,
\chi_k) \right|
\]
  and let
$$F_{T}(x)=\frac{1}{T}\textup{meas}\left\{\tau\in[0, T] : g(\tau)<x\right\}$$
be a distribution function of $ g(\tau)$. If
$g(\tau)$ is $B.A.P.$ then it is known (see
Jessen and Wintner \cite[Theorem 27]{jw} or
Laurin\v cikas \cite[Theorem 6.3, Chapter
2]{laur}) that there is a distribution function
$F(x)$ such that $F_{T}(x)$ converges weakly to
$F(x)$ for $T\to \infty$. It means that if $F(x)$
is continuous at $x=\eps$ then
$$\lim_{T\to\infty} F_T(\eps)$$
exists. Thus to obtain Theorem \ref{th2} we need
to show that $g(\tau)$ is $B.A.P.$

We remark that if $a(t)$ and $b(t)$ are both
non-negative $B.A.P.$ functions of $t$, then,
$t\mapsto \max (a(t),b(t))$ is also $B.A.P.$
since $\max (a(t),b(t))$ can be written as
\begin{align*}
\max (a(t),b(t))=\frac{1}{2}\left( \left|
a(t)-b(t)\right| +(a(t)+b(t))\right) ,
\end{align*}
and the modulus of $B.A.P.$ function is again
$B.A.P.$ By this we have that  $g(\tau)$ is
$B.A.P.$ if the function
$$f(\tau)=\max_{s\in\mathcal
K}\left|L(s+id_1\tau, \chi_1)-L(s+id_2\tau,
\chi_2)\right|$$ is $B.A.P.$ In view of the note
below the formula (\ref{bap}) the function
$f(\tau)$ is $B.A.P.$ if there are
 $U.A.P$ functions $f_N(\tau)$ such that
\begin{align}\label{uap}
 \lim_{N \to +\infty }\left( \lim
\sup_{T\rightarrow +\infty
}\frac{1}{2T}\int_{-T}^{T}\left|
f(\tau)-f_N(\tau)\right|^{2}d\tau\right) =0.
\end{align}

Let
$$L_N(s,\chi)=\sum _{n\le N}\frac{\chi(n)}
{n^s}$$ be a partial sum of the Dirichlet series
associated with $L(s,\chi)$. Next we show that
the equality (\ref{uap}) is true with
\[f_N(\tau)=\max_{s\in\mathcal
K}\left|L_N(s+id_1\tau,\chi_1)- L_N(s+id_2\tau,
\chi_2)\right|.\]

 By repeating the proof of
Proposition 12 of Mauclaire \cite{mauc} we get
that $f_N(\tau)$ is $U.A.P.$ for any
$d_1,d_2\in\Bbb R$. Note that the case when $d_1$
or $d_2$ is equal to zero is already included in
Proposition 12 of Mauclaire \cite{mauc}.

Further we have that
\begin{align*}
&L(s+id_1\tau,\chi_1)-L(s+id_2\tau,\chi_2)\\
&=\left(L(s+id_1\tau,\chi_1)-L_N(s+id_1\tau,\chi_1)+L_N(s+id_2\tau,\chi_2)-L(s+id_2\tau,\chi_2)\right)\\
&\phantom{=}+\left(L_N(s+id_1\tau,\chi_1)-L_N(s+id_2\tau,\chi_2)\right),
\end{align*}
and as a consequence, we get that
\begin{align*}
&\left| f(\tau) -f_N(\tau) \right| \\
&\leq
\sup_{s\in K}\left|
L(s+id_1\tau,\chi_1)-L_N(s+id_1\tau,\chi_1)+L_N(s+id_2\tau,\chi_2)-L(s+id_2\tau,\chi_2)\right|
\\ &\le \sup_{s\in K}\left|
L(s+id_1\tau,\chi_1)-L_N(s+id_1\tau,\chi_1)\right|\\
&\phantom{\le}+\sup_{s\in
K}\left|L_N(s+id_2\tau,\chi_2)-L(s+id_2\tau,\chi_2)\right|.
\end{align*}
Then, in view of the inequality
$(a+b)^2\le2a^2+2b^2$, we obtain that
\begin{align*}
&\frac{1}{2T}\int_{-T}^{T}\left| f(\tau) -f_N(\tau) \right| ^{2}dt \\
&\leq \frac{1}{T}\int_{-T}^T\left( \sup_{s\in
K}\left| L(s+id_1\tau,\chi_1)-L_N(s+id_1\tau,\chi_1)\right| \right)
^{2}dt\\&\phantom{\le}+\frac{1}{T}\int_{-T}^T\left(
\sup_{s\in K}\left|
L(s+id_2\tau,\chi_2)-L_N(s+id_2\tau,\chi_2)\right| \right) ^{2}dt.
\end{align*}
By Mauclaire \cite[Theorem 5.1]{mauc2} we have
that, for any real $d$,
\[
\lim_{N \to +\infty }\left( \lim
\sup_{T\rightarrow +\infty }\frac{1}{2T^{\prime
}}\int_{-T}^{T}\left( \sup_{s\in K}\left|
f(s+idt)-f_N(s+idt)\right| \right) ^{2}dt\right)
=0.
\]
This proves the equality (\ref{uap}) and Theorem \ref{th2}
\endproof

From the proof we see that Theorem \ref{th2}
remains true with Dirichlet $L$-functions
$L(s,\chi_j)$, $j=1,\dots,m$,  replaced by any
general Dirichlet series satisfying conditions of
Theorem 5.1 of Mauclaire \cite{mauc2}.

{\bf Acknowledgment.} We thank Jean-Loup
Mauclaire for suggesting Theorem~\ref{th2} and
for other useful comments which helped to improve
the paper.

{\bf Remark.} The `$\liminf$' version of
Corollary \ref{cor} is  independently obtained by
Takashi Nakamura in ``The generalized strong
recurrence for non-zero rational parameters",
arXiv:1006.1778v1 [math.NT].


\begin{thebibliography}{99}

\bibitem{bagdis}{\sc B. Bagchi}, {\it The statistical behaviour and universality properties of the Riemann zeta-
function and other allied Dirichlet series}, PhD
Thesis, Calcutta: Indian Statistical Institute,
1981.

\bibitem{bag}{\sc B. Bagchi}, {\it A joint universality theorem for Dirichlet L-functions},
Math. Z., {\bf 181} (1982), 319-334.

\bibitem{bag87}{\sc B. Bagchi}, {\it Recurrence in topological dynamics and the Riemann hypothesis},
Acta Math. Hung., {\bf 50} (1987), 227-240.

\bibitem{bak} {\sc A. Baker}, {\it Transcendental number theory}, London: Cambridge University Press. X,
(1975).

\bibitem{bes} {\sc A.S. Besicovitch}, {\it Almost periodic functions},   Dover, New York, (1954).

\bibitem{bohr}{\sc H. Bohr}, {\it \"Uber eine quasi-periodische Eigenschaft Dirichletscher Reihen mit
    Anwendung auf die Dirichletschen $L$-Funktionen},
Math. Ann., {\bf 85} (1922), 115-122.

\bibitem{jw} {\sc B. Jessen and A. Wintner},
{\it Distribution functions and the Riemann zeta
function},   Trans. Am. Math. Soc. 38 (1935),
48-88.

\bibitem{gs}{\sc E. Girondo and J. Steuding}, {\it Effective estimates for the distribution of values of Euler products},
Monatsh. Math., {\bf 145}, No. 2 (2005), 97-106.

\bibitem{gon} {\sc S.M. Gonek}, {\it Analytic properties
of zeta and L-functions}, Ph. D. Thesis,
University of Michigan, 1979.

\bibitem{kls}{\sc J. Kaczorowski, A. Laurin\v cikas,  and J. Steuding}, {\it On the value distribution of shifts of universal Dirichlet series},
Monatsh. Math., {\bf 147} (2006), 309-317.

\bibitem{vor} {\sc A.A. Karatsuba and S.M. Voronin}, {\it The Riemann
zeta-function}, De Gruyter Expositions in
Mathematics. 5. Berlin etc.: W. de Gruyter. xii,
(1992).

\bibitem{laur}{ \sc A. Laurin\v cikas}, {\it Limit theorems for the Riemann zeta-function},
 Mathematics and its Applications (Dordrecht). 352. Dordrecht: Kluwer Academic Publishers, (1995).

\bibitem{mauc}{\sc J.-L. Mauclaire}, {\it Almost periodicity and Dirichlet series}, Laurin\v cikas, A. (ed.) et al., Analytic and probabilistic methods in number theory. Proceedings of the 4th international conference in honour of J. Kubilius, Palanga, Lithuania, September 25--29, 2006. Vilnius: TEV, (2007),109-142.

\bibitem{mauc2}{\sc J.-L. Mauclaire}, {\it On some Dirichlet series},
 Proceedings of the conference “New Directions in the Theory of Universal Zeta- and L-Functions”, W\" urzburg, Germany, October 6-10, 2008. Shaker Verlag, (2009), 171-248.

\bibitem{nak}{\sc T. Nakamura}, {\it The joint universality and the generalized strong recurrence for Dirichlet $L$-functions},
Acta Arith., {\bf 138} (2009), 357-362.


\bibitem{pan}{\sc \L. Pa\'nkowski}, {\it Some remarks on the generalized strong recurrence for $L$-functions},
in: New Directions in Value Distribution Theory
of zeta and L-Functions: proceedings of Wurzburg
Conference, October 6-10, 2008, Shaker Verlag,
(2009), 305-315.

\bibitem{steu}{\sc J. Steuding}, {\it Value distribution of
$L$-functions}, Lecture Notes in Mathematics
1877, Springer, 2007.

\bibitem{titch}{\sc E.C. Titchmarsh},
{\it The theory of functions. 2nd ed.}, London:
Oxford University Press. X,  1975.

\end{thebibliography}
\end{document}